\documentclass[12pt,a4paper]{article}
\usepackage{amsmath,amscd,amssymb,amsfonts,amsthm}
\usepackage[mathscr]{eucal}
\usepackage{hyperref}

\def\co{\colon\thinspace}
\newcommand{\mathcall}{\EuScript}

\renewcommand{\leq}{\leqslant}

\newtheorem{thm}{Theorem}[section]

\newtheorem{lm}{Lemma}[section]
\newtheorem{prop}{Proposition}[section]
\newtheorem{cor}{Corollary}[section]

\theoremstyle{definition}
\newtheorem{de}{Definition}[section]
\newtheorem{ex}{Example}[section]

\newtheorem{rem}{Remark}[section]

\DeclareMathOperator{\Ber}{Ber} 
 \DeclareMathOperator{\ord}{ord}
\DeclareMathOperator{\grad}{grad} \DeclareMathOperator{\Res}{Res}
\renewcommand{\div}{\mathop{\mathrm{div}}}
\newcommand{\divg}{\mathop{\mathrm{div}}\nolimits_{\gamma}}
\newcommand{\divgg}{\mathop{\mathrm{div}}\nolimits_{\bar\gamma}}
\newcommand{\divv}[2]{\div\nolimits_{\frac{\g_{#1}+\g_{#2}}{2}}{(\g_{#2}-\g_{#1})}}
\newcommand{\coc}[2]{c(\g_{#1},\g_{#2})}

\newcommand{\Schw}{{\mathfrak{S}}}
\newcommand{\BVAct}{{\mathscr{S}}}
\newcommand{\Act}{{\mathscr{A}}}

\newcommand{\gn}{{\gamma}_{\#}}

\newcommand{\bro}{{\boldsymbol{\rho}}}

\DeclareMathOperator{\Vol}{Vol} 
 
 \DeclareMathOperator{\Der}{Der}
\DeclareMathOperator{\Vect}{Vect} \DeclareMathOperator{\sub}{sub}
\DeclareMathOperator{\symb}{symb} \DeclareMathOperator{\ad}{ad}
\newcommand{\lie}[1]{{\cal L}_{{#1}}}
\newcommand{\liew}[1]{{\cal L}^w_{{#1}}}

\newcommand{\volx}{{Dx}}
\newcommand{\der}[2]{{\frac{\partial {#1}}{\partial {#2}}}}

\newcommand{\RR}{\mathbb R}
\newcommand{\Z}{{\mathbb Z_{2}}}

\newcommand{\p}{\partial}

\newcommand{\fun}{C^{\infty}}
\newcommand{\V}{{\mathfrak{V}}}

\def\e{\varepsilon}
\def\s{\sigma}
\def\f{{\varphi}}

\def\D{\Delta}
\def\G {\Gamma}
\def\o{\omega}
\newcommand{\g}{{\gamma}}
\newcommand{\F}{{\Phi}}
\renewcommand{\S}{{{S}}}
\newcommand{\lt}{\theta} % last term in the master Hamiltonian for long bracket

\renewcommand{\l}{{\lambda}}
\newcommand{\x}{{\xi}}
\def\d{\delta}
\def\t{\theta}

\newcommand{\df}{{\mathcall{D}\f}}
\newcommand{\dff}{{\mathcall{D}\F}}

\newcommand{\ft}{{\tilde f}}
\newcommand{\kt}{{\tilde k}}

\newcommand{\at}{{\tilde a}}
\newcommand{\bt}{{\tilde b}}
\newcommand{\ct}{{\tilde c}}

\newcommand{\Xt}{{\tilde X}}
\newcommand{\At}{{\tilde A}}
\newcommand{\Bt}{{\tilde B}}
\newcommand{\Yt}{{\tilde Y}}

\newcommand{\bs}{{\boldsymbol{s}}}

\newcommand{\psit}{{\tilde \psi}}
\newcommand{\chit}{{\tilde \chi}}

\newcommand{\ps}{{\boldsymbol{\psi}}}
\newcommand{\ch}{{\boldsymbol{\chi}}}
\newcommand{\ph}{{\boldsymbol{\varphi}}}

\newcommand{\newga}{{\bar\g^a}}
\newcommand{\newlt}{{\bar\lt}}

\newcommand{\X}{{\boldsymbol{X}}}

\title{On odd Laplace operators. II}

\newcommand{\myaddress}{{%\vspace{-0.4cm}
\footnotesize $^{\scriptstyle{1}}$ Department of Mathematics,
University of Manchester Institute of Science and Technology (UMIST),
PO Box 88, Manchester M60 1QD, England\\ \smallskip} {\footnotesize
$^{\scriptstyle{2}}$ G.~S.~Sahakian~Department~of~Theoretical
~Physics, Yerevan State University, \\1 A. Manoukian Street, 375049
Yerevan, Armenia\\ \smallskip} {\footnotesize $^{\scriptstyle{3}}$
Joint Institute for Nuclear Research, Dubna 141980, Russia\\}
{\footnotesize {\tt theodore.voronov@umist.ac.uk,
khudian@umist.ac.uk} }}

\author{Hovhannes M. Khudaverdian{\small $^{\scriptstyle{1,2,3}}$},
Theodore  Voronov{\small $^{\scriptstyle{\,1}}$}}
\date{\myaddress}

\begin{document}
\maketitle
\vspace{-1.0cm}
\begin{abstract}\footnotesize
\noindent We analyze geometry of the second order differential
operators, having in mind applications to Batalin--Vilkovisky
formalism in quantum field theory.  As we show, an exhaustive picture
can be obtained by considering pencils of differential operators
acting on densities of all weights simultaneously. The algebra of
densities, which we introduce here, has a natural invariant scalar
product. Using it, we prove that there is a one-to-one correspondence
between second-order operators in this algebra and the corresponding
brackets. A bracket on densities incorporates a bracket on functions,
an ``upper connection'' in the bundle of volume forms, and a term
similar to the ``Brans--Dicke field'' of the Kaluza--Klein formalism.
These results are valid for even operators on a usual manifold as
well as for odd operators on a supermanifold. For an odd operator
$\Delta$ we show that conditions on the order of the operator
$\Delta^2$ give an hierarchy of properties such as flatness of the
upper connection and the Batalin--Vilkovisky master equation. In
particular, we obtain a complete description of generating operators
for an arbitrary odd Poisson bracket.
\end{abstract}

\tableofcontents

\section{Introduction}

This paper is a direct continuation of our paper~\cite{tv:laplace1}.
However, it is completely independent and can be read as such.

There are two motivations for this work.

The first motivation comes from the Batalin--Vilkovisky formalism in
quantum physics. The problem is to give a description of all the
so-called  `$\Delta$-operators'. A particular algebraic problem
related to this is to give a description of all generating operators
for a given odd Poisson algebra. The second motivation is entirely
geometrical. The problem is to describe geometric structures encoded
in a differential operator. In particular, the question is, what is
necessary to recover a differential operator from its principal
symbol.

It is well known that various quantities in quantum field theory can
be expressed via the Feynman integral
\begin{equation}   \label{eqfeynint1}
    Z=\int e^{\frac{i}{\hbar}\,S[\f]}\, \df.
\end{equation}
Here integration is over all field configurations and $S[\f]$ stands
for the classical action functional.   However, if the theory
possesses a gauge freedom (like electrodynamics or Yang--Mills
theory), a  modification is required. According to the most
up-to-date comprehensive procedure,
--- the Batalin--Vilkovisky formalism, --- the recipe is as
follows~\cite{bv:perv, bv:vtor, bv:closure}. The
(infinite-dimensional) manifold of fields $\f$ has to be extended and
necessarily becomes  a supermanifold possessing an odd symplectic
structure. The classical action $S[\f]$ and the
integral~\eqref{eqfeynint1} are  replaced by an ``extended'' action
$\BVAct[\F]$ and by  the integral\footnote{Strictly speaking, there
must be $\sqrt{\dff}$ instead of ${\dff}$, since it is the
half-densities that give volume forms on Lagrangian
submanifolds~\cite{hov:max, hov:proclms}.}
\begin{equation}   \label{eqfeynint2}
    Z_{\text{BV}}=\int e^{\frac{i}{\hbar}\,\BVAct[\F]}\, \dff\,,
\end{equation}
where integration is over a Lagrangian submanifold. The key condition
is that the extended action must satisfy the Batalin--Vilkovisky
``quantum master equation''
\begin{equation}   \label{eqquantumbv}
   \D e^{\frac{i}{\hbar}\,\BVAct[\F]} =0,
\end{equation}
and this secures the gauge invariance of the quantum theory. The
Batalin--Vilkovisky $\D$-operator is  an ``odd Laplacian'' associated
with an odd symplectic structure. We want to emphasize that the
precise geometric formulation of the above-said, including the
precise definition of the operator $\D$, is a non-trivial task, to
which a lot of work has been devoted, see~\cite{hov:delta, hov:khn1,
hov:khn2, hov:max, hov:semi, hov:proclms, tv:laplace1} and
\cite{ass:bv, ass:symmetry, ass:semiclassical}. Initially $\D$ has
been thought of as an operator acting on functions. In~\cite{hov:max,
hov:semi} it was shown that $\D$ should be considered as an operator
on half-densities (= semidensities) rather than functions.  Moreover,
as we have shown in~\cite{tv:laplace1}, the best understanding of
$\D$ can be achieved by considering it on densities of various
weights. There is a remarkable similarity between odd Poisson
geometry and the usual Riemannian geometry noticed
in~\cite{tv:laplace1}.

The Batalin--Vilkovisky operator ``generates'' (in a precise
sense) the odd symplectic structure of the extended phase space.
In a more abstract setup, every `$\D$-operator' on functions
generates an odd bracket. A question that remained open, is how to
describe all operators generating a given bracket. A lot of work
was devoted to this problem, in the geometric as well as in an
algebraic setting. (Notice that the Batalin--Vilkovisky
quantization formalism motivated the introduction of various
algebraic structures, such as the ``Batalin--Vilkovisky
algebras'', the study of which has developed into an independent
area.) See in particular~\cite{bt:1993gen, bt:1993multi},
\cite{hov:khn1}, \cite{yvette:divergence}.  The current paper
gives a complete solution of this problem.

A related  question is how to obtain an operator on half-densities or
densities of any other weight from an operator on functions, or vice
versa.

These questions naturally bring us to the second, purely geometrical
motivation for the present paper.

Suppose we are given a differential operator $\D$ of order $\leq n$
acting on functions on a manifold $M$. Which geometric  structures
are naturally associated with it? We want to stress that we are
considering an operator acting on functions on a manifold without any
extra structure (like a Riemannian structure) given \textit{a
priori}. If we write this operator in local coordinates as
\begin{equation}\label{eqdifop}
    \D=
    \sum_{k=0}^n
    \frac{1}{k!}\,A^{a_1,\ldots,a_k}(x)\,\p_{a_1}\ldots\p_{a_k},
\end{equation}
the coefficients $A^{a_1,\ldots,a_k}$ are transformed in a
complicated way under a change of coordinates. Which geometric
information is encoded in them?

First of all, as it is well known, the top order part defines an
invariant object, called the principal symbol of $\D$:
\begin{equation}\label{eqprincsymb}
    \s (\D)=\frac{1}{n!}\,
    A^{a_1,\ldots,a_n}(x)\,p_{a_1}\ldots p_{a_n}.
\end{equation}
The principal symbol $\s (\D)$ is an invariantly defined function on
$T^*M$. It also can be approached   from a purely algebraic viewpoint
(see subsection~\ref{subsechigh}).

What about the geometric meaning of the lower order terms
in~\eqref{eqdifop}? H\"ormander has introduced the so-called
subprincipal symbol of $\D$, which  essentially is the principal
symbol of $\D-(-1)^n\D^*$ where the adjoint operator $\D^*$
(acting on the same space as $\D$)  depends on a choice of volume
element, e.g., given by a coordinate system. Hence the
subprincipal symbol $\sub \D$ is not a genuine function on $T^*M$,
but has a non-trivial transformation law. In the standard approach
this is considered as a nuisance that can be overcome by trading
functions for half-densities, for which the adjoint operator  is
intrinsically defined. However, a closer look at the
transformation law of the subprincipal symbol $\sub \D$ for
operators on functions reveals that it is very similar to a
connection in the bundle of volume forms $\Vol M$.

This is quite unexpected as, let us repeat, we have started from an
operator acting   on scalar functions, with no geometric data (like
bundles and connections) being given beforehand. In particular, this
prompts to consider operators acting on densities of various weights
$w\in \RR$ together. As soon as we adopt this viewpoint, the picture
immediately clears.

We introduce the \textit{algebra of densities} $\V(M)$ on a
manifold $M$ as the algebra of densities of all weights $w\in\RR$
under the tensor multiplication. Formally, it is the algebra of
sections of the direct sum of the bundles $|\Vol (M)|^w$ over all
$w\in \RR$. $\V(M)$ possesses a natural invariant scalar product,
and it contains, in particular, functions, volume forms and
half-densities. One can consider differential operators in the
algebra $\V(M)$, and there is a natural notion of the adjoint
operator. It is possible to give a nice classification of
derivations of $\V(M)$.

The main results of the present paper can be summarized as follows.

We describe differential operators in $\V(M)$. A differential
operator of order $\leq 2$ in $\V(M)$ is equivalent to a quadratic
pencil
\begin{equation*}
    \D_w=\D_0+w\,A+w^2\,B,
\end{equation*}
where $\D_0$ is a second-order operator on functions, $A$ and $B$
have order $1$ and $0$ respectively.  We prove that there is a
one-to-one correspondence between the self-adjoint operators with the
condition $\D_0(1)$=0 and the corresponding brackets in $\V(M)$. A
bracket in $\V(M)$ from the viewpoint of $M$ is a ``long bracket''
incorporating a bracket on functions, an ``upper connection'' in the
bundle of volume forms and a term similar to the ``Brans--Dicke
field'' $g^{55}$ of the Kaluza--Klein formalism in the relativity
theory. This   gives a complete description of the geometric
information necessary to recover an operator from the corresponding
bracket on functions (= the principal symbol) and unfolds the
relations between operators acting on densities of different weights.

These results are valid for even operators   as well as for odd
operators on a supermanifold. For odd operators and brackets,
analysis can be pursued   further, leading to remarkable results
that have no analogs in the even case. For an odd operator
$\Delta$ it is natural to study the operator
$\D^2=\frac{1}{2}\,[\D,\D]$. We show that   conditions on the
order of  $\Delta^2$ give an hierarchy of properties, including
the flatness condition for the upper connection and the
Batalin--Vilkovisky master equation.

In particular, this gives a complete description of the generating
operators  for an arbitrary odd Poisson bracket.

Though we mainly consider operators of order $\leq 2$ here, we
discuss how our results can be generalized to differential operators
of higher order. A generalization is also possible to operators of
non-zero weight $\l$ (i.e., those mapping densities of weight $w$ to
densities of weight $w+\l$).

The structure of the paper is the following.

In Section~\ref{secoperbrack} we study arbitrary operators of the
second order on functions and densities, and the corresponding
brackets. We interpret the subprincipal symbol as an upper
connection. We define the algebra $\V(M)$ and establish its
properties. We prove the main classification theorem giving a $1-1$
correspondence between operators and brackets in $\V(M)$, and
consider examples.

In Section~\ref{secjacobi} we consider odd operators and odd
brackets. We study the Jacobi identity for an odd  long bracket  and
conditions on the operator $\D^2$ corresponding to various properties
of brackets.

In Section~\ref{secgen} we discuss generalizations.

\section{Second orders operators and long  brackets} \label{secoperbrack}

\subsection{Subprincipal symbol as a connection}\label{subsecsubprinc}

Let $\D$ be a differential operator acting on functions on a manifold
$M$. Its principal symbol is a well-defined function on the cotangent
bundle $T^*M$. Compared to it, the so-called subprincipal symbol
$\sub \D$ introduced by H\"or\-man\-der is not a genuine function on
$T^*M$ but depends on a choice of coordinates. (Only for operators
acting on half-densities the subprincipal symbol becomes an invariant
function.) We want to point out that  the transformation law for the
subprincipal symbol $\sub \D$ for an operator on functions or on
densities of any weight $w\neq \frac{1}{2}$ allows to interpret it as
a sort of connection.

Indeed, let $\D$ be an operator of the second order. (Operators of
higher order will be treated later.) To avoid complications with
signs let $\D$ be even and $M$ be an ordinary manifold. Suppose in
local coordinates
\begin{equation}\label{eqlap}
    \D=\frac{1}{2}S^{ab}\,\p_b\p_a+ T^a\,\p_a +R,
\end{equation}
where $S^{ab}=S^{ba}$. Then (with standard conventions) the principal
symbol of $\D$ is $\symb \D=-S^{ab}p_bp_a$, which is an invariant
quadratic function on $T^*M$, and the subprincipal symbol in local
coordinates is~\cite{hoer:volume3}:
\begin{equation}\label{eqsub}
    \sub \D=i\left(2T^a-\p_bS^{ba}\right)p_a.
\end{equation}
It is nothing but the principal symbol of $\D-\D^*$ where $\D^*$ is a
coordinate-dependent adjoint operator, defined using a coordinate
volume element $Dx$.  Denote
\begin{equation}\label{eqsub2}
    \g^a=\p_bS^{ba}-2T^a,
\end{equation}
so that $\sub \D=-i\g^ap_a$. (In the sequel we omit the factors $i$.)
Then under a change of coordinates the coefficients $\g^a$ are
transformed as follows:
\begin{equation}\label{eqzamenagamma}
    \g^{a'}=\left(\g^a+S^{ab}\,\p_b\log J\right)\der{x^{a'}}{x^a},
\end{equation}
where $J=Dx'/Dx$ is the Jacobian. If we assume that the matrix
$S^{ab}$ is invertible, so that $S_{ab}S^{bc}=\d_a{}^c$, then it
is possible to lower indices and get $\g_a=S_{ab} \g^b$ with the
transformation law
\begin{equation}\label{eqzamenagamma2}
    \g_{a'}=\left(\g_a+\p_a\log J\right)\der{x^{a}}{x^{a'}}.
\end{equation}
This is the transformation law for the coefficients of a
connection in the bundle of volume forms on $M$. We see that
\textit{the subprincipal symbol of an operator of the second order
acting on functions defines a connection in $\Vol M$}, with
\begin{equation}\label{eqconn}
    \nabla_a\rho=(\p_a +\g_a)\rho,
\end{equation}
$\g_a=S_{ab}(\p_cS^{cb}-2T^b)$, if the ``upper metric'' $S^{ab}$
given by the principal symbol is non-degenerate. In general, the
subprincipal symbol of $\D$ defines a so-called \textit{upper
connection} or ``contravariant derivative'' in the bundle $\Vol
M$:
\begin{equation}\label{equpconn}
    \nabla^a\rho=(S^{ab}\p_b +\g^a)\rho,
\end{equation}
over the map $S^{\#}\co T^*M\to TM$ defined by the principal
symbol.

``Upper connections'' or ``contravariant derivatives'' were
considered earlier, in particular, in the context of Poisson
geometry. A natural framework for them is  that of Lie algebroids
(vector bundles with a Lie bracket of sections and a Lie
homomorphism of sections to vector fields). However, for our
purposes this framework is not entirely suitable, since in general
we do not expect to have a Lie bracket of covector fields,  ---
except for the case of odd Poisson geometry, see later. Instead,
we shall use as an appropriate formalism the following language of
``long brackets''. (``Long'' is meant to remind ``long
derivatives'', the physicists' term for covariant derivatives.)

Fix a bi-derivation in $\fun(M)$, denoted as a bracket $\{\ ,\
\}$.  We do not assume  the Jacobi identity. Let $E\to M$ be a
vector bundle. A \textit{long bracket} between functions and
sections of $E$ over a given bracket of functions is a bilinear
operation $(f, \bs)\mapsto \{f;\bs\}$, where $f\in \fun(M)$ is a
function and $\bs\in \fun(M, E)$ is a section, taking values in
sections of $E$, with the properties:
\begin{align}
    \{fg;\bs\}&=f\{g;\bs\}+\{f;\bs\}g, \label{eqlong1}\\
    \{f;g\bs\}&=\{f,g\}\bs+g\{f;\bs\}.
\end{align}
Equation~\eqref{eqlong1} means that the value of $\{f;\bs\}$
depends only on  $df$. A long bracket is related with an ``upper
connection'' by the formula $\nabla^{df}\bs:=\{f;\bs\}$.

With an operator on functions $\D$ given by formula~\eqref{eqlap}
we can associate the following  bracket of functions  and a  long
bracket  over it:
\begin{align}
    \{f,g\}&:=S^{ab}\p_b f\,\p_a g  \label{eqcoorbracket}\\
    \{f;\bro\}&:=\left(S^{ab}\p_b f\,\p_a \rho +\g^b\p_b \label{eqcoorlongbracket}
    f\, \rho\right)\volx,
\end{align}
where $\bro=\rho \volx$ is a volume form. Here $\g^a$ are given
by~\eqref{eqsub2}. (If $S^{ab}$ is a metric, then the bracket
$\{f,g\}$ is the scalar product of gradients.) An alternative
coordinate-free  expression is as follows:
\begin{align}
    \{f,g\}&:=\D(fg)-\D f\,g-f\,\D g +\D(1)\,fg, \label{eqinvbracket} \\
    \{f;\bro\}&:=\D(f\bro)-\D f\,\bro-f\,\D \bro
    +\D(1)\,f\bro,\label{eqinvlongbracket}
\end{align}
where in~\eqref{eqinvlongbracket} we define $\D$ on volume forms
as the adjoined operator $\D^*$. A direct check shows
that~(\ref{eqinvbracket},\ref{eqinvlongbracket})
give~\eqref{eqcoorbracket},\eqref{eqcoorlongbracket}
and~\eqref{eqsub2}; in particular,  the coordinate-free
formulae~(\ref{eqinvbracket},\ref{eqinvlongbracket}) yield a proof
of the transformation law~\eqref{eqzamenagamma}.

Let us summarize.

{Any second order differential operator $\D$ acting on {functions}
defines a  bracket  of functions~\eqref{eqinvbracket} ---
essentially, the ``polarized'' principal symbol of $\D$ --- and  an
``upper connection''~\eqref{eqsub2} in the bundle of {volume forms}
$\Vol M$ --- essentially, the subprincipal symbol of $\D$
--- which may be written as a ``long
bracket''~\eqref{eqinvlongbracket} between functions and
volume forms extending the bracket of functions.}

Hence, starting from operators acting on functions, we are
naturally prompted to consider densities. Moreover, the  long
bracket  $\{f;\bro\}$ defined above can be extended by a ``Leibniz
rule'' from densities of weight $1$ (volume forms) to densities of
arbitrary weight. For our purposes we shall need a further
generalization with the first argument also replaced by a density.
As we shall see below, this will give a ``completion'' of the
theory: after extending both operators and brackets to arbitrary
densities it will become possible to establish a one-to-one
correspondence between them.

\subsection{The algebra of densities}

Let $M$ be a supermanifold. Consider densities of arbitrary weights
$w\in \RR$. (To avoid complications, we can assume that all
appropriate Berezinians between local charts are positive.) Under the
tensor product  densities form a graded commutative algebra:
$\ps\ch=(-1)^{\psit\chit}\ch\ps$, $w(\ps\ch)=w(\ps)+w(\ch)$. Tilde
stands for parity, $w$ for weight; we drop $\otimes$ from the
notation. Denote the algebra of densities on $M$ by $\V(M)$.

The commutative algebra $\V(M)$ can, in fact, be identified with a
certain algebra of functions on an extended  manifold $\hat M:=(\Ber
TM)\setminus M$, i.e., the frame bundle of the determinant bundle
$\Ber TM$. The natural coordinates on $\hat M$ induced by local
coordinates $x^a$ on $M$ are $x^a, t$ where $t$ can be identified
with the volume element $Dx$. A formal sum of densities of various
weights $\ps=\sum \psi_w(x)(Dx)^w\in\V(M)$ can be identified with its
``generating function'' $\ps(x,t)=\sum \psi_w(x)t^w$ (summation over
a finite number of weights). In the sequel we shall use elements of
$\V(M)$ and the corresponding functions on $\hat M$ interchangeably.

The algebra of densities $\V(M)$ has a natural bilinear scalar
product:
\begin{equation}\label{eqscprod}
    \langle\ps,\ch\rangle:=\int_M \psi\chi Dx
\end{equation}
if $w(\ps)+w(\ch)=1$; otherwise the scalar product is zero. In terms
of functions on $\hat M$ the scalar product can be expressed as
\begin{equation}\label{eqscprod2}
    \langle\ps,\ch\rangle:=\int_M
    \Res\left(t^{-2}\,\ps(x,t)\,\ch(x,t)\right)
    Dx,
\end{equation}
where $\Res$ stands for the residue at $t=0$. Notice that this
algebra possesses a unit. The scalar product satisfies the invariance
condition
\begin{equation*}%\label{eqscprod}
    \langle\ps\ch,\ph\rangle=\langle\ps,\ch\ph\rangle.
\end{equation*}
One can consider (formally) adjoint operators w.r.t. the scalar
product~\eqref{eqscprod}. In particular, one has $t^*=t$,
$\p_t^*=-\p_t+2t^{-1}$, $\p_a^*=-\p_a$; all functions of $x^a$ are
self-adjoint. We have $\hat w^*=1-\hat w$.

Let us describe $\Der\V(M)$, i.e.  vector fields on $\hat M$. A
vector field on $\hat M$ of weight $\l$ has the form
\begin{equation}\label{eqvectoronhatm}
    \X = t^{\l}\left(X^a(x)\,\der{}{x^a}+X_0(x)\,t\der{}{t}\right),
\end{equation}
or
\begin{equation}\label{eqvectoronhatm2}
    \X = t^{\l}\bigl(X^a(x)\,\p_a+X_0(x)\,\hat w\bigr).
\end{equation}
The first term has the meaning of a vector density of weight $\l$ on
$M$, i.e., a derivation taking functions to densities of weight $\l$
(this is just the restriction of $\X$ to functions on $M$). The
second term has no independent meaning if $X^a\neq 0$;  notice the
transformation law
\begin{align*}
    X^{a'}&=J^{-\l}X^a\,\der{x^{a'}}{x^a}\\
    X_0'&=J^{-\l}\left(X_0+X^a\p_a\log J\right).
\end{align*}

There is a canonical operation $\div$ (the divergence) on vector
fields on $\hat M$. This is no wonder, since there is an invariant
scalar product of functions on $\hat M$, in other words a generalized
volume form. The explicit formula for the divergence is
\begin{equation}\label{eqdivfordensities}
    \div
    \X=t^{\l}\Bigl(\p_aX^a(-1)^{\at(\Xt+1)}+(\l-1)\,X_0\Bigr),
\end{equation}
if $\X$ is given by~\eqref{eqvectoronhatm}, \eqref{eqvectoronhatm2}.

\begin{thm}\label{thmclassderonhatm}
For $\l\neq 1$, every derivation of weight $\l$ in the algebra of
densities $\V(M)$ can be uniquely decomposed into the sum of a
divergence-free derivation and a derivation of the form
$t^{\l}f(x)\hat w$; every divergence-free derivation has the form
\begin{equation}\label{eqdivfreederivation}
    \X=t^{\l}\left(X^a(x)\,\p_a-\frac{1}{\l-1}\,\p_aX^a(-1)^{\at(\Xt+1)}\,\hat w\right)
\end{equation}
and is uniquely defined by a vector density $X=(Dx)^{\l}X^a\,\p_a$ on
$M$. The decomposition of a general $\X$ is
\begin{equation}\label{eqdecompositionofder}
    \X=\underbrace{t^{\l}\left(X^a(x)\,\p_a-
    \frac{1}{\l-1}\,\p_aX^a(-1)^{\at(\Xt+1)}\,\hat w\right)}_{\text{divergence-free
    part}}+ \frac{1}{\l-1}\,\div\X\,\hat w.
\end{equation}
\end{thm}

For $\l=0$ the derivation~\eqref{eqdivfreederivation} is   the Lie
derivative along a vector field $X=X^a(x)\,\p_a$ on $M$ acting on
densities. In general, the divergence-free
derivation~\eqref{eqdivfreederivation} can be seen as a
``generalized Lie derivative'' along the vector density
$X=(Dx)^{\l}X^a\,\p_a$.

Part of what we are doing is of a purely algebraic nature and holds
for arbitrary commutative associative algebras with a unit and with
an invariant (non-degenerate) scalar product. In such an algebra, if
an operator $\D$ is of  order $k$, in the algebraic sense, then the
adjoint operator $\D^*$ is also of order $k$; the difference
$\D-(-1)^k\D^*$ is of order $k-1$.   For derivations of such an
algebra there exists a canonical divergence.

Recall that in any commutative associative algebra $A$ an abstract
divergence operator  is a linear map $\div\co \Der A\to A$ with the
property
\begin{equation}\label{eqmaindiv}
    \div (aX)=a\div X+(-1)^{\at\Xt} X(a)
\end{equation}
(see~\cite{koszul:crochet85}, \cite{yvette:divergence}). For
functions on a manifold $M$  an abstract divergence and a
connection in $\Vol M$ are equivalent notions. Indeed, such an
equivalence is established by the formula
\begin{equation}\label{eqdivandconnection}
    \int_M  (\div X) \, \bro =- \int_M  \nabla_X\bro,
\end{equation}
and the property~\eqref{eqmaindiv} translates into the characteristic
property  of a covariant derivative. (Of course, this can be
transported into a more abstract setting.) The explicit formulae:
\begin{equation*}%\label{}
    \nabla_a\rho=(\p_a+\g_a)\,\rho, \quad \
    \div X=(\p_a-\g_a)X^a (-1)^{\at(\Xt+1)}.
\end{equation*}
For abstract divergence operators there is a notion of ``curvature''
(see~\cite{yvette:divergence}). This is exactly the curvature of the
corresponding connection on volume forms.

\begin{rem}\label{remlapwithgenuineconnect}
Second order differential operators on functions for which the
associated upper connections  come from genuine connections as
$\g^a=S^{ab}\g_b$ are exactly the ``Laplacians''
\begin{equation*}%\label{}
    \D f=\frac{1}{2}\,(\p_a-\g_a)\left(S^{ab}\p_bf\right)=
    \frac{1}{2}\,\divg\grad f,
\end{equation*}
where $\grad f=S^{ab}\p_b f$ and $\div_{\g}$ is the divergence
operator corresponding to the connection $\g_a$.
\end{rem}

If in the algebra $A$ there is an invariant scalar product, we can
define a canonical operation $\div$ by either of the equivalent
formulae: for $X\in \Der A$,
\begin{equation}\label{eqdefdiv1}
    \langle\div X, a\rangle=-\langle 1, X(a)\rangle,
\end{equation}
or
\begin{equation}\label{eqdefdiv2}
    \div X=-(X+X^*).
\end{equation}
Here $X^*$ stands for the  operator formally adjoint to $X$.  Notice
that  $X+X^*$ is of order $0$, i.e., an element of $A$. Immediately
checked is that  $\div$  is an even, linear operation
satisfying~\eqref{eqmaindiv}. Also,  an identity
\begin{equation}\label{eqdivkommutator}
    \div [X,Y]=X(\div Y)+(-1)^{\Xt\Yt} Y(\div X),
\end{equation}
holds, meaning that the curvature  of the operator $\div$ is zero.
\begin{ex}
For  derivations of $\V(M)$ we easily get from~\eqref{eqdefdiv1},
\eqref{eqdefdiv2} the explicit formula~\eqref{eqdivfordensities}.
\end{ex}

\begin{rem}
Formula~\eqref{eqdefdiv2} is very close to the usual definition of a
subprincipal symbol.  Define $\div$ acting on arbitrary operators of
order $\leq k$ and taking them to operators of order $\leq k-1$ as
$\div \D:=-(\D-(-1)^k\D^*)$. Then it easily follows that
\begin{equation*}
    \div [\D_1,\D_2]=[\div \D_1, \D_2]+[\D_1,\div \D_2] +[\div \D_1,\div \D_2],
\end{equation*}
which implies~\eqref{eqdivkommutator} for derivations.  In this
abstract setting, the  {subprincipal symbol} of $\D$ can be
defined as the principal symbol of $\div \D$, i.e., as the class
of $\div \D$ modulo operators of order $\leq k-2$. Notice also
that $\div^2=0$, so we are getting a complex.
\end{rem}

\subsection{Equivalence between operators and brackets}\label{subsecequiv}

In this section we shall prove that any bracket on $\hat M$, i.e.,
a bracket in the algebra $\V(M)$, is in a $1-1$ correspondence
with a differential operator of the second order in $\V(M)$. This
should be compared with the fact that a bracket of functions gives
only the principal symbol of the generating operator  and does not
allow to recover it in full.

In any algebra (commutative associative with a unit) by a
\textit{bracket} we mean a symmetric bi-derivation:
\begin{align}
    \{a,b\}&=(-1)^{\at\bt}\{b,a\}  \label{eqsym}\\
    \{ka,b\}&=(-1)^{\kt\e}k \{a,b\} \label{eqlin}\\
    \{a,bc\}&=\{a,b\}c+(-1)^{(\at+\e)\bt}b\{a,c\}. \label{eqleibn}
\end{align}
Here  $\e\in \Z$ is the parity of the bracket. Let us emphasize
that in the context of this paper we consider  symmetric brackets
and hence exclude the usual (antisymmetric even) Poisson brackets.

Notice that from~(\ref{eqsym}-\ref{eqleibn}) follows
$\{ab,c\}=(-1)^{\at\e}a\{b,c\}+(-1)^{\bt\ct}\{a,c\}b$. In our
notation the parity $\e$ ``sits'' at the opening bracket.

\begin{de}
A  \textit{long bracket} on $M$ is a bracket in the algebra $\V(M)$.
Notation: $\{\ps;\ch\}$.
\end{de}

We use semi-colon for a long bracket, reserving comma for a bracket
of functions. A long bracket has weight $\l\in \RR$ if
$w(\{\ps;\ch\})=w(\ps)+w(\ch)+\l$. For $\l=0$ we have a
generalization of ``long brackets'' $\{f;\bs\}$ considered in
section~\ref{subsecsubprinc}.

Since a long bracket on $M$ is a usual bracket from the viewpoint of
$\hat M$, it can be specified by a master Hamiltonian
$\hat\S\in\fun(T^*\hat M)$ as
\begin{equation}\label{eqabstractlongbracket}
    \{\ps;\ch\}=((\hat\S,\ps),\ch)
\end{equation}
(see~\cite{tv:laplace1}). The parentheses denote the canonical
Poisson bracket on the cotangent bundle.  For a bracket of parity
$\e$ and weight $\l$ the master Hamiltonian has  to be of the form
\begin{equation}\label{eqmasterham}
    \hat \S=t^{\l}\frac{1}{2}\left(\S^{ab}p_bp_a+ 2t\g^ap_ap_t+t^2\lt
    p_t^2\right),
\end{equation}
where $p_a$ and $p_t$ are the momenta conjugate to $x^a$ and $t$
respectively.  $\hat S$ is of parity $\e$. That means that for a
long bracket we have the following expression:
\begin{multline}\label{eqexplicittlongbracket}
    \{\ps;\ch\}=t^{\l}\left(\S^{ab}\der{\ps}{x^b}\der{\ch}{x^a}(-1)^{\at\psit} \right.
    \\
    \left.+t\g^a\left(\der{\ps}{x^a}\der{\ch}{t}+(-1)^{\at\psit}\der{\ps}{t}\der{\ch}{x^a}\right)
     +t^2\lt \der{\ps}{t}\der{\ch}{t}\right).
\end{multline}
Notice that $\hat w=t\p_t$ is the weight operator taking eigenvalues
$w$ on densities of weight $w$.
Formula~\eqref{eqexplicittlongbracket} can be rewritten using $\hat
w$ applied to $\ps$ and $\ch$.  Taking off the hats, we come to an
equivalent description of the long
bracket~\eqref{eqexplicittlongbracket} as a ``double pencil'' of
brackets $\{\ \,;\ \}_{w_1,w_2}$:
\begin{multline}\label{eqlongbracketpencil}
    \{\psi (Dx)^{w_1};\chi (Dx)^{w_2}\}_{w_1,w_2}=
    \left(\S^{ab}\,\p_b{\psi}\p_a{\chi}(-1)^{\at\psit} \right.
    \\
    \left.+\g^a\left(w_2\,\p_a{\psi}\,{\chi}+(-1)^{\at\psit}w_1\,{\psi}\,\p_a{\chi}\right)
     + w_1w_2\,\lt\,{\psi}{\chi}\right)(Dx)^{w_1+w_2+\l}.
\end{multline}
We shall often suppress the subscripts $w_1,w_2$ in the notation.
Equivalently,
\begin{align*}
    \{x^a;x^b\}&=\S^{ab}(Dx)^{\l}\\
    \{x^a;Dx\}&=\g^a(Dx)^{\l+1}\\
    \{Dx;Dx\}&=\lt (Dx)^{\l+2}.
\end{align*}
It is useful to have the formulae for the transformation law of the
coefficients of a long bracket under a change of coordinates. One can
deduce that
\begin{align}
    \S^{a'b'}&=J^{-\l}\,\S^{ab}\,\der{x^{b'}}{x^b}\der{x^{a'}}{x^a}\,(-1)^{\bt'(\at'+\at)}
    \label{eqzamenaSl}
    \\
    \g^{a'}&=J^{-\l}\left(\g^a+\S^{ab}\,\p_b\log
    J\right)\der{x^{a'}}{x^a}
    \label{eqzamenagammal}
    \\
    \lt'&=J^{-\l}\left(\lt+2\g^a\p_a\log J +\S^{ab}\,\p_b\log J\,\p_a\log J\right)
    \label{eqzamenathetal}
\end{align}
where $J=Dx'/Dx$ is the Jacobian (Berezinian). We shall be mainly
interested in the case $\l=0$. Then it follows that a long bracket
incorporates a bracket of functions as well as an ``upper
connection'' $\g^a$. The space of all long brackets is a vector
space. For a fixed bracket of functions given by $S^{ab}$, upper
connections $\g^a$ form an affine space  associated with the vector
space  of vector fields on $M$. Similarly, for $S^{ab}$, $\g^a$
fixed, the coefficients $\lt$ make up an affine space associated with
the space of functions on $M$. See examples later.
%The meaning of the coefficient $\lt$ will be elucidated later.

\begin{rem}
If we write the components of the tensor on $\hat M$ specifying a
long bracket as a block matrix, then
\begin{equation}\label{eqkaluzamatrix}
    (\hat S^{\hat a\hat b})=
    \begin{pmatrix}
       t^{\l}S^{ab} & t^{\l+1} \g^a   \\
      t^{\l+1} \g^a & t^{\l+2}\lt   \\
    \end{pmatrix}
\end{equation}
and we can see a straightforward analogy with the Kaluza--Klein
formalism in field theory, where a metric tensor in a $5$-dimensional
spacetime combines the usual metric tensor together with a gauge
field and an extra scalar field (the ``Brans--Dicke field'' of
tensor-scalar theories of gravitation).
\end{rem}

Now we are going to formulate the central theorem of this paper. It
has an abstract algebraic counterpart, which is almost trivial if
properly stated.

Let $A$ be a commutative associative algebra with a unit. Consider an
operator $\D$ in $A$ of parity $\e$ and introduce an operation
$\{a,b\}$ by the formula
\begin{equation}\label{eqmainbracket}
     \{a,b\}=\D(ab)-\D a\,b-(-1)^{\at\e}a\,\D b +\D(1)\,ab.
\end{equation}
Clearly, this operation has parity $\e$, is bilinear and symmetric:
$\{a,b\}=(-1)^{\at\bt}\{b,a\}$, $\{ka,b\}=(-1)^{\kt\e}k\{a,b\}$. We
shall call $\D$ a {\em generating operator}
for~\eqref{eqmainbracket}.

\begin{prop}
\label{thmmainabs} $\vphantom{hernia}$

{\em (1) }   The operation~\eqref{eqmainbracket} is a bracket, i.e.,
satisfies~\eqref{eqleibn},  if and only if $\D$ is of order $\leq 2$,
in the algebraic sense. Two operators generate the same bracket if
and only if they differ by an operator of the first order.

{\em (2) } Suppose that the algebra $A$ possesses an invariant scalar
product. Given a bracket in $A$, a generating operator for it is
uniquely determined by the extra conditions $\D^*=\D$ and $\D(1)=0$.
\end{prop}

The first statement is essentially due to
Koszul~\cite{koszul:crochet85}. Recall  that $\ord \D\leq 2$ if
all the triple commutators $[[[\D,a],b],c]$ vanish. The derivation
property of~\eqref{eqmainbracket} w.r.t. each argument turns out
to be equivalent to this condition.  As for the second statement,
a proof is straightforward: let $\D'$ be an arbitrary operator
generating a given bracket, then one can check that
$\D:=\D''-\D''(1)$, where $\D'':=\frac{1}{2}(\D'+{\D'}^*)$, is a
generating operator and it is the only generating operator
satisfying the conditions $\D^*=\D$ and $\D(1)=0$.

For the algebra of densities $\V(M)$  we can find the generating
operator explicitly.

What is an  operator of the second order in the algebra $\V(M)$? It
is convenient to use the language of the extended manifold $\hat M$
and then translate back using the weight operator $\hat w$. Every
operator in $\V(M)$ is equivalent to a pencil of operators $\D_w$
acting on $w$-densities (and mapping them to $(w+\l)$-densities for
operators/pencils of weight $\l$).

\begin{lm}\label{lemsecondorder}
An operator of order $\leq 2$ in the algebra $\V(M)$ is equivalent to
a quadratic pencil of the form
\begin{equation}\label{eqpencil}
    \D_w=\D_0+wA+w^2B,
\end{equation}
where $\D_0$ is an operator of order $\leq 2$ acting on functions,
$A$ and $B$ are operators of order $\leq 1$ and $0$ respectively.
\textnormal{(Note that $A$, $B$ do not make sense independently of
$\D_0$.)}
\end{lm}

In the language of pencils, the \textit{adjoint pencil} for any
pencil $\D_w$, corresponding to the adjoint of the operator in
$\V(M)$, is given by $(\D_{1-\l-w})^*$, because $\hat w^*=1-\hat
w$.

\begin{thm}\label{thmmain}
For a given long bracket on $M$, there exists a unique operator
$\D$ on $\V(M)$ that satisfies
\begin{equation}\label{eqcond}
    \D^*=\D, \quad \D(1)=0,
\end{equation}
i.e., $\D_w^*=\D_{1-w}$, $\D_0(1)=0$,  and generates the
bracket by the formula
\begin{equation}\label{eqmainbracket1}
     \{\ps;\ch\}=\D(\ps \ch)-\D \ps\cdot\ch-(-1)^{\psit\e}\ps\cdot\D
     \ch.
\end{equation}
If the bracket is given by~\eqref{eqexplicittlongbracket}
or~\eqref{eqlongbracketpencil}, then $\D$ is given by
\begin{multline}\label{eqdeltaforbracket}
    \D=t^{\l}\frac{1}{2}\Biggl(\S^{ab}\p_b\p_a+
    \left(\p_b\S^{ba}(-1)^{\bt(\e+1)}+(2\hat w+\l-1)\g^a\right)\p_a +\\
    \hat w\, \p_a\g^a(-1)^{\at(\e+1)}  +
        \hat w(\hat w+\l-1)\,\lt \Biggr).
\end{multline}
\end{thm}

Rewriting formula~\eqref{eqdeltaforbracket} as an operator pencil,
we get:
\begin{multline}\label{eqdeltaforbracket2}
    \D_w=\frac{1}{2}\Biggl(\S^{ab}\p_b\p_a+
    \left(\p_b\S^{ba}(-1)^{\bt(\e+1)}+(2 w+\l-1)\g^a\right)\p_a +\\
     w\,\p_a\g^a(-1)^{\at(\e+1)}  +
         w(w+ \l-1)\,\lt \Biggr).
\end{multline}
We call~\eqref{eqdeltaforbracket2} the \textit{canonical pencil}
corresponding to the long bracket~(\ref{eqexplicittlongbracket},
\ref{eqlongbracketpencil}).

Consider three important ``values'' of this pencil (where we set
$\l=0$). The operator acting on functions ($w=0$)
\begin{equation}\label{eqdeltaonfunct}
    \D_0=\frac{1}{2}\left(\S^{ab}\p_b\p_a+\left(\p_b\S^{ba}(-1)^{\bt(\e+1)}-\g^a\right)\p_a\right);
\end{equation}
the operator acting on volume forms ($w=1$)
\begin{equation}\label{eqdeltaonvol}
    \D_1=\frac{1}{2}\left(\S^{ab}\p_b\p_a+\left(\p_b\S^{ba}(-1)^{\bt(\e+1)}+\g^a\right)\p_a
    +\p_a\g^a(-1)^{\at(\e+1)}\right);
\end{equation}
the operator acting on half-densities ($w=\frac{1}{2}$)
\begin{equation}\label{eqdeltaonhalf}
    \D_{1/2}=\frac{1}{2}\left(\S^{ab}\p_b\p_a+\p_b\S^{ba}(-1)^{\bt(\e+1)}\p_a +
    \frac{1}{2}\,\p_a\g^a(-1)^{\at(\e+1)} -\frac{1}{4}\,\lt \right).
\end{equation}
All of them have the same principal symbol defined by the bracket
on functions.  Notice that operators acting on functions and on
volume forms do not depend on the coefficient $\lt$. Given the
principal symbol, they are completely defined by the ``upper
connection'' $\g^a$. The operator on half-densities, on the other
hand, depends on $\g^a$ only via $\p_a\g^a(-1)^{\at(\e+1)}$;
instead, it includes $\lt$ in its zeroth-order term. Neither
operator allows to recover the long bracket by itself; knowing the
two $\D_0$ and $\D_{1/2}$ or $\D_1$ and $\D_{1/2}$ is sufficient.
Outside of these exceptional cases it is sufficient to know any
single operator $\D_{w_0}$ to recover $\g^a$ and $\lt$, i.e. the
whole pencil $\D_w$ (see Theorem~\ref{thmpencilbyindop} below).

Thus a self-adjoint operator  of the second order (in the
algebraic sense) on densities, vanishing on $1$, contains and is
completely defined by the following data: a bracket on functions,
an ``upper connection'' $\g^a$ on volume forms, and a quantity
$\lt$ analogous to the Brans--Dicke field. We shall elucidate the
meaning of this quantity in the examples below.

Theorem~\ref{thmmain} can be proved by a direct calculation
starting from a general expression for a pencil $\D_w$ with
indeterminate coefficients depending on $w$. The
property~\eqref{eqpencil} will come about automatically.

There is an alternative approach showing  that the generating
operator for a given long bracket is nothing but the Laplace operator
``$\div\grad$''. Though this  adds nothing in terms of explicit
formulae, it will be useful for deducing further properties of $\D$.

Let us for a moment come back  to a general algebraic situation. Let
$A$ be an algebra with a scalar product as above. Consider  a bracket
in $A$. Define $\grad a:=\{a,\ \}$. This is a derivation of $A$ of
parity $\at+\e$. The operator $\grad\co A\to \Der A$ is linear of
parity $\e$. We have
\begin{equation}\label{eqgradproduct}
    \grad (ab)= (-1)^{\e\at}a \,(\grad
    b) +(-1)^{(\at+\e)\bt}b\,(\grad a).
\end{equation}
We can introduce the \textit{Laplace operator} in the algebra $A$,
corresponding to a given bracket, as
\begin{equation}\label{eqdefabstractlaplace}
    \D a:=\div \grad a.
\end{equation}
It is an operator of parity $\e$.

\begin{prop}
The Laplace operator satisfies $\D(1)=0$, $\D^*=\D$, and
\begin{equation}\label{eqlaplaceasgenerator}
    \D (ab):=\D a\,b+(-1)^{\e\at}a\,\D b +2\{a,b\}
\end{equation}
for all $a$, $b$.
\end{prop}
\begin{proof}
Notice that $\grad (1)=0$. By~\eqref{eqdefdiv1}, we have
$\langle\D a,b\rangle=-\langle 1, \{a,b\}\rangle$, which implies
$\langle\D a,b\rangle=(-1)^{\at\e}\langle a,\D b\rangle$ due to
the symmetry of the bracket. Identity~\eqref{eqlaplaceasgenerator}
follows by applying~\eqref{eqmaindiv} to~\eqref{eqgradproduct}.
\end{proof}

It follows that up to a factor of $2$,~\eqref{eqdefabstractlaplace}
is   the unique generating operator for the bracket given by
Proposition~\ref{thmmainabs}. Thus, for algebras with an invariant
scalar product we have the uniqueness and existence theorem for
generating operators:

\begin{thm}\label{thmlaplacian}
Let $A$ be a commutative associative   algebra with unit and an
invariant scalar product. For a given bracket in $A$, the unique
generating operator which is self-adjoint and vanishes on constants
is, up to a factor of $2$, the Laplace operator $\div\grad$.
\end{thm}

In particular, in the algebra $\V(M)$ we have for $\ps=t^{w}\psi(x)$
\begin{equation}\label{eqgraddens}
    \grad \ps=t^{w+\l}(-1)^{\at\psit}\bigl(S^{ab}\p_b\ps +\g^a
    w\,\psi\bigr)\,\p_a +t^{w+\l+1}\bigl(\g^a\p_a\ps +\lt w
    \psi\bigr)\,\p_t.
\end{equation}
Using the formula for  the divergence~\eqref{eqdivfordensities}, we
after simplification get
\begin{multline*}
    \D\ps=\frac{1}{2}\,\div\grad\ps=\frac{t^{\l}}{2}\,\Bigl(S^{ab}\p_b\p_a+
    \left(\p_bS^{ba}(-1)^{\bt(\e+1)}+(2w+\l-1)\g^a\right)\p_a+ %\right.
    \\
    %\left.
    w\p_a\g^a(-1)^{\t(\e+1)}+w(w+\l-1)\lt\Bigr)(t^{w}\psi)
\end{multline*}
(where we have restored the factor $1/2$), and finally we arrive at
formulae~\eqref{eqdeltaforbracket}, \eqref{eqdeltaforbracket2}.

\subsection{Properties of the canonical pencil. Examples}

Let us study the change of the canonical pencil under a change of
$\g^a$, $\lt$ with $S^{ab}$ fixed. Suppose $\newga=\g^a+X^a$,
$\newlt=\lt+\x$. Let  $\bar\D_w$ denote  the pencil corresponding to
$\newga$, $\newlt$, and $\D_w$ denote  the pencil corresponding to
$\g^a$, $\lt$. Let $\D$ and $\bar\D$ be the corresponding operators
in the algebra $\V(M)$. Since the canonical pencil depends linearly
on the data $S^{ab}, \g^a, \lt$, the difference $\bar\D_w-\D_w$ is
the canonical pencil corresponding to the long bracket given by the
matrix
\begin{equation*}%\label{}
    \begin{pmatrix}
      \text{\large $\boldsymbol{0}$}& t^{\l+1}X^a \\
      t^{\l+1}X^a & t^{\l+2}\x \\
    \end{pmatrix}.
\end{equation*}
Immediately follows that $X^a$ and $\x$ transform as
\begin{align*}
    X^{a'}&=J^{-\l} X^a \der{x^{a'}}{x^a}\\
    \x'&=J^{-\l}\left(\x+2X^a\p_a\log J\right).
\end{align*}
Hence
\begin{equation*}%\label{}
    \X=X^a\p_a+\frac{1}{2}\,\x\,\hat w
\end{equation*}
is a vector field on $\hat M$. Recall that there is a canonical
divergence on $\hat M$:
\begin{equation}\label{eqdivxbold}
    \div\X=\p_a X^a(-1)^{\at(\e+1)}+(\l-1)\frac{1}{2}\,\x.
\end{equation}

\begin{thm}\label{thmzamenagammaiteta}
Under a change of $\g^a$, $\lt$ with $S^{ab}$ fixed we have
\begin{equation}\label{eqnewdelta}
    \bar\D-\D=\frac{1}{2}\,\X\,(2\hat w+\l-1)+\frac{1}{2}\,\div\X\,\hat
    w,
\end{equation}
or, in terms of pencils:
\begin{equation}\label{eqnewdelta2}
    \bar\D_w-\D_w=\frac{1}{2}\,(2w+\l-1)\left(X^a\p_a+\frac{1}{2}\,w\x\right)+
    \frac{1}{2}\,w\,\div\X.
\end{equation}
\end{thm}

\begin{cor}
Decomposing $\X$ for $\l\neq 1$ we obtain:
\begin{equation}\label{eqnewdelta3}
    \bar\D_w-\D_w=\left(w-\frac{\l+1}{2}\right)\liew{X}+
    \frac{w(w-1)}{\l-1}\,\div\X
\end{equation}
where by
\begin{equation}\label{eqliew}
    \liew{X}:=X^a\p_a-\frac{w}{\l-1}\,\p_aX^a(-1)^{\at(\e+1)}
\end{equation}
we denoted the divergence-free part of $\X$.
\end{cor}

In the case $\l=0$ which is particularly interesting for us,
\eqref{eqliew}  is the usual Lie derivative $\lie{X}$ of
$w$-densities along the vector field $X$ on $M$, and we have
\begin{equation}\label{eqnewdelta4}
    \bar\D_w-\D_w=\left(w-\frac{1}{2}\right)\lie{X}-
     {w(w-1)}\,\div\X
\end{equation}
where the canonical divergence of $\X=X^a\p_a+(1/2)\x\hat w$ is
expressed as
\begin{equation*}%\label{}
    \div\X=\p_aX^a(-1)^{\at(\e+1)}-\frac{1}{2}\,\x.
\end{equation*}

In the rest of this subsection we work with operators of weight
$\l=0$.

Consider important constructions of canonical pencils.

As follows from~\eqref{eqnewdelta4}, the subspace of canonical
pencils with   $S^{ab}=0$ is the direct sum of two natural subspaces:
$V_1=\{(2w-1)\lie{X} \ \,|\, X\in \Vect (M)\}$ and $V_2=\{w(w-1)f  \
\,|\, f\in \fun(M)\}$. The following example provides a subspace
which is complementary to $V_1\oplus V_2$.

\begin{ex}[Canonical pencil associated with a volume form]\label{exLB}
\label{exdeltarho} Fix a basis volume form $\bro=\rho Dx$.  For a
given bracket of functions specified  by a tensor $\S^{ab}$, the
following ``Laplace--Beltrami type'' formula
\begin{equation}\label{eqbeltrami}
    \D^{LB}(\psi
    (Dx)^w):=\frac{1}{2}\rho^{w-1}\p_a\bigl(\rho\S^{ab}\p_b(\rho^{-w}
    \psi)\bigr)(-1)^{\at(\e+1)}\,(Dx)^w
\end{equation}
defines a self-adjoint operator $\D^{LB}$ on densities, vanishing on
the unit. Expanding~\eqref{eqbeltrami} and comparing
with~\eqref{eqdeltaforbracket2}, it is easy to see that it gives a
canonical pencil with the following $\g^a$ and $\lt$:
\begin{align*}
    \g^a&=-\S^{ab}\p_b\log \rho\\
    \lt&= \S^{ab}\p_b\log \rho\cdot \p_a\log \rho=\g^a\g_a.
\end{align*}
Here the upper connection comes from a genuine connection
$\g_a:=-\p_a\log\rho$ in $\Vol M$ (which is flat).
\end{ex}

It was exactly the pencil $\D^{LB}_w$ of Example~\ref{exdeltarho}
that was the main object of study in~\cite{tv:laplace1}.

Using the canonical pencil $\D^{LB}_{w}$ associated with a volume
form it is possible to give a convenient parametrization of
\textit{all} canonical pencils. If we fix a volume form $\bro$, then
every canonical pencil has the appearance
\begin{equation}\label{eqcanoncherezlaplace}
   \D_w=\D^{LB}_w+\frac{1}{2}\,(2w-1)\lie{X}+w(w-1)f
\end{equation}
for some vector field $X$  and a scalar function $f$. Conversely,
every pencil~\eqref{eqcanoncherezlaplace} is a canonical pencil. This
is a decomposition   (depending on a choice of $\bro$) of the space
of canonical pencils into the direct sum of three subspaces. If
$\D_w$ is given by~\eqref{eqcanoncherezlaplace}, then
\begin{align*}
    \g^a&=-S^{ab}\p_b\log\rho+X^a\\
    \lt&=S^{ab}\p_b\log\rho\,\p_a\log\rho + 2\p_aX^a(-1)^{\at(\e+1)}+f.
\end{align*}

We can apply this to obtain  a coordinate-free decomposition of an
arbitrary second order linear differential operator acting on
densities of fixed weight $w_0$ on a manifold $M$.

Any such operator in local coordinates has the appearance
\begin{equation}\label{eqoperatoronwdens}
    L=\frac{1}{2}\,S^{ab}\,\p_b\p_a+T^a\p_a+R.
\end{equation}
Clearly, the principal symbol $S^{ab}$ defines a bracket on $M$; from
the subprincipal symbol of $L$ one can construct an upper connection
$\g^a$ in $\Vol M$ by the formula
\begin{equation}\label{eqgammapooperatoru}
    \g^a=\frac{1}{2w_0-1}\,\left(2T^a-\p_bS^{ba}(-1)^{\bt(\e+1)}\right) \quad
    \text{if $w_0\neq \frac{1}{2}$}.
\end{equation}
In the case  $w_0=1/2$ the subprincipal symbol
$2T^a-\p_bS^{ba}(-1)^{\bt(\e+1)}$ is a vector field, so it does not
give any upper connection.

If we choose a volume form $\bro=\rho Dx$, then we can write $L$ as
\begin{equation}\label{eqlcherezdelta}
    L=\D^{LB}_{w_0}+\lie{Q}+f
\end{equation}
where $Q$ is a vector field and  $f$ is a scalar function on $M$,
both uniquely defined by $L$. $\lie{Q}$ stands for a Lie derivative
on $w_0$-densities. Indeed, let $\G_a=-\p_a\log \rho$ be the
connection in $\Vol M$ generated by $\bro$, and $\G^a=S^{ab}\G_b$ be
the corresponding upper connection. The vector field $Q$ is defined
by the formula
\begin{equation}\label{eqQ}
    Q^a=\frac{1}{2}\,(2w_0-1)\left(\g^a-\G^a\right)=
    \frac{1}{2}\,\left(2T^a-\p_bS^{ba}(-1)^{\bt(\e+1)}-
    (2w_0-1)\G^a\right).
\end{equation}
(Notice that $Q$ makes sense for all $w_0$ including $w_0=1/2$.) Then
one can see that $L-\D^{LB}_{w_0}-\lie{Q}$ is an operator of the
zeroth order, i.e. the multiplication by a scalar function, hence we
get~\eqref{eqlcherezdelta}.  Formula~\eqref{eqlcherezdelta} can be
viewed as a replacement of the coordinate
description~\eqref{eqoperatoronwdens} where instead of coefficients
depending on a coordinate system we consider a vector field $Q$ and a
scalar function $f$ depending on a volume form $\bro$.

What is the relation of an individual operator $L$ with canonical
pencils?

We can use decompositions~\eqref{eqcanoncherezlaplace} and
\eqref{eqlcherezdelta}. The general picture is as follows. There
is a specialization map $\D_w\mapsto \D_{w_0}=L$ from the linear
space of all canonical pencils $\D_w$ to the linear space of all
second order differential operators on densities of a fixed weight
$w_0$. This map is an isomorphism for all non-singular values
$w_0\neq 0,1, \frac{1}{2}$\,:

\begin{thm}\label{thmpencilbyindop}
Let $L$ be a second order differential operator acting on densities
of weight $w_0$. If $w_0\neq 0,1, \frac{1}{2}$, then there exists a
unique canonical pencil $\D_w$ passing through $L$, i.e.
$L=\D_{w_0}$. Namely,
\begin{equation}\label{eqpencilcherezl}
    \D_w=\D^{LB}_w + \frac{2w-1}{2w_0-1}\,\lie{Q}+
     \frac{w(w-1)}{w_0(w_0-1)}\,f
\end{equation}
where $Q$ and $f$ are given by~\eqref{eqlcherezdelta}.
\end{thm}

For $w_0=\frac{1}{2}$ the image of the specialization map consists of
all self-adjoint operators (recall that the space of half-densities
has a scalar product) and the kernel is the subspace
$V_1=\{(2w-1)\lie{X}\}$. For $w_0=0$ the image of the specialization
map consists of all operators vanishing on constants and the kernel
is the subspace $V_2=\{w(w-1)f\}$.  Similarly we can describe the
specialization map for $w_0=1$.

\begin{ex}[Canonical pencil associated with a connection]
\label{exdeltasviaznost} The construction of Example~\ref{exdeltarho}
can be generalized for an arbitrary connection $\g_a$ in $\Vol M$:
given a bracket of functions, both $\g^a$ and $\lt$ can be defined by
the connection as $\g^a=S^{ab}\g_b$ and $\lt=\g^a\g_a$. In a more
abstract language, the operator $\D_w$ is (up to $1/2$) the Laplace
operator $\div \grad$ on $w$-densities, where $\grad$ is the
covariant gradient w.r.t. the induced connection $w\g_a$, and
$\div=\div_{\g}$ is the divergence of vector fields on $M$ defined by
the connection $\g_a$ in $\Vol M$. In particular, if a linear
connection in $TM$ is given, there is an associated connection in the
bundle $\Vol M$, namely $\g_a=-\G_{ab}^b(-1)^{\bt}$, in the standard
notation, with the curvature given by the trace of the Riemann tensor
$-R_{abc}^c(-1)^{\ct}$. (Notice that this gives an example of a
divergence operator with a possibly non-zero curvature.)
\end{ex}

Examples~\ref{exdeltarho} and~\ref{exdeltasviaznost} explain the
geometrical meaning of $\lt$.

Consider the restriction of the general formulae for the
transformation of the canonical pencil under a change of $\g^a$,
$\lt$ to canonical pencils defined by a connection $\g_a$. A change
of connection is given by a covector field $X_a$:
\begin{equation}\label{eqchangelowergamma}
    \bar\g_a=\g_a+X_a.
\end{equation}
Since $\lt$ is defined by $\g_a$ as $\g^a\g_a$, we get for its change
\begin{equation}\label{eqchangethetabylowergamma}
    \bar\lt=\lt+2\g^aX_a+X^aX_a.
\end{equation}

\begin{prop} For operators defined by a bracket on $M$ and a connection
$\g_a$, the change of the operator under a change of connection is
given by the formula
\begin{equation}\label{eqchangeofdeltaunderconnection}
    \bar \D-\D=\frac{1}{2}\,(2w-1)\lie{X}-
    w(w-1)\left(\divg X-\frac{1}{2}\,X^a X_a\right)
\end{equation}
where $\bar\g_a=\g_a+X_a$, $X^a=S^{ab}X_b$.
\end{prop}

It follows that covector fields $X$ considered as functions of two
connections: $X_a=X_a(\g,\bar\g)=\bar\g_a-\g_a$, possess the
following groupoid property. If for three connections $\g,\bar \g,
\bar\bar\g$ with $\bar\g_a=\g_a+X_a$, $\bar\bar\g_a=\bar\g_a+Y_a$ the
equations
\begin{align}
    \divg X-\frac{1}{2}\,X^a X_a&=0 \label{eqbv1}\\
    \divgg Y-\frac{1}{2}\,Y^a Y_a&=0  \label{eqbv2}
\end{align}
are satisfied, then the equation
\begin{equation}
    \divg (X+Y)-\frac{1}{2}\,(X+Y)^a (X+Y)_a=0  \label{eqbv3}
\end{equation}
is satisfied for $X+Y$. Equations~\eqref{eqbv1}, \eqref{eqbv2},
\eqref{eqbv3} are generalization  of the ``Batalin--Vilkovisky
equations'' of the ``master groupoid'' discovered
in~\cite{tv:laplace1}. It follows that the specialization of $\D_w$
to $w_0=\frac{1}{2}$ (operators on half-densities) does not depend on
a connection $\g_a$ but only on its orbit w.r.t. a groupoid action:
$\g_a\mapsto \g_a+X_a$ where $X$ satisfies~\eqref{eqbv1}.

\begin{rem}\label{remtwistedvolumeform}
A special feature of the connection in Example~\ref{exdeltarho} is
flatness. Every flat connection in $\Vol M$ is locally represented by
a $1$-form $\g_a=-\p_a\log\rho$. On the intersections we get
$\log\rho -\log\rho'=\log J +c$, where $c$ is a local constant.
Clearly, $c$ is a $1$-cocycle. The local functions $\rho$ can be
glued to a nonvanishing volume form if and only if the cohomology
class $[c]\in H^1(M;\RR)$ equals zero. We may say that every flat
connection in $\Vol M$ comes from a  volume form up to the described
``twist''.
\end{rem}

It is tempting to relate flatness of a connection with properties of
the operator $\D$ and the long bracket. This we shall do in
Section~\ref{secjacobi}.

\section{Jacobi identity and flatness}\label{secjacobi}

In this section we shall focus on the case of odd operators and odd
brackets. For them the general analysis performed above can be
advanced further. The main tool of classification will become the
operator $\D^2$. There is a sharp contrast with even operators and
the corresponding  even brackets, for which no similar development is
possible.

\subsection{Jacobi identity for a long bracket}
In the previous sections, a ``bracket''   meant just a bilinear
concomitant satisfying the derivation property w.r.t. each of its
arguments. Now we want to move further and explore the possibility of
imposing a Jacobi identity. It turns out, however, that this will be
possible only for {odd} brackets. Indeed, notice that all brackets
considered above are symmetric
--- compared to Lie brackets, which are antisymmetric.
What kind of a Jacobi identity can be introduced for a symmetric
bracket? The usual Jacobi identity in Lie algebras can be
reformulated as either of the following properties: the linear map
$a\mapsto [a,\ ]$ takes the bracket to the commutator of
operators, or the operator $[a, \ ]$ for each element of the
algebra is a derivation of the bracket. The symmetry and linearity
conditions for our brackets have the form
\begin{align}%
    \{a,b\}&=(-1)^{\at\bt}\{b,a\} \label{eqsymbrack} \\
    \{ka,b\}&=(-1)^{\kt\e}k\{a,b\}, \  \{a,bk\}=\{a,b\}k;
\end{align}
hence it would make sense to consider a modified bracket
$[a,b]:=(-1)^{\at\e}\{a,b\}$, for which the same conditions read
as
\begin{align*}%\label{}
    [a,b]&=(-1)^{\at\bt+\at\e+\bt\e}[b,a]\\
    [ka,b]&=k[a,b], \ [a,bk]=[a,b]k.
\end{align*}
Now, for an even bracket we have $[a,b]=\{a,b\}$, and  symmetry:
$[a,b]=(-1)^{\at\bt}[b,a]$,  still holds. Since the commutator of
operators is antisymmetric, $[A,B]=-(-1)^{\At\Bt}[B,A]$, there is no
hope for $a\mapsto \ad a=[a,\ ]$ to take brackets to brackets. Is it
possible, however, to have $\ad a$ as a derivation of the even
bracket? Suppose this is satisfied:
\begin{equation}\label{eqfakejacobi}
    [a,[b,c]]=[[a,b],c]+(-1)^{\at\bt}[b,[a,c]]
\end{equation}
for all $a,b,c$. Rearranging cyclically and adding with suitable
signs, we arrive at the following

\begin{prop}
If for an even symmetric bracket the ``fake Jacobi''
property~\eqref{eqfakejacobi} is satisfied, then all triple brackets
vanish: $[[a,b],c]=0$.
\end{prop}

In the differential-geometric situation for a bracket
$[f,g]=S^{ab}\p_bf\p_ag$ with a symmetric even tensor $S^{ab}$, this
implies $S^{ab}\equiv 0$.

Conclusion: there is no way of imposing Jacobi identity for an even
symmetric bracket. Hence, we have to concentrate on odd brackets.

\textbf{From now on all brackets are odd.}

For an odd symmetric bracket $\{a,b\}$, we have
$[a,b]=(-1)^{\at}\{a,b\}$, and the symmetry
condition~\eqref{eqsymbrack} for $\{a,b\}$ becomes the antisymmetry
condition for $[a,b]$ with a shift of parity:
\begin{equation*}
[a,b]=-(-1)^{(\at+1)(\bt+1)}[b,a]
\end{equation*}
The standard Jacobi condition w.r.t. the shifted parity for
$[a,b]=(-1)^{\at}\{a,b\}$,  translates  into the condition
\begin{equation}\label{eqjacobiasder}
    \{a,\{b,c\}\}=(-1)^{\at+1}\{\{a,b\},c\}+
    (-1)^{(\at+1)(\bt+1)}\{b,\{a,c\}\},
\end{equation}
or
\begin{equation}\label{eqjacobicycle}
    (-1)^{\at\ct}\{\{a,b\},c\}+
    (-1)^{\ct\bt}\{\{c,a\},b\}+(-1)^{\bt\at}\{\{b,c\},a\}
    =0
\end{equation}
for the symmetric bracket $\{a,b\}$.  Notice that in $\{a,b\}$ the
left opening bracket is odd. Historically, when Lie superalgebras
first appeared in topology, they were written using a symmetric
bracket, and the Jacobi identity  appeared for them exactly in the
form~\eqref{eqjacobicycle}.

In the sequel we  continue to work with the symmetric brackets.
(Notice the different choice of signs in~\cite{tv:laplace1}.)

After this algebraic digression, let us return to our geometric
situation. Suppose on a manifold $M$ is given an odd symmetric
bracket. It is specified by the master Hamiltonian as
\begin{equation}\label{eqoddbracket}
    \{f,g\}=((S,f),g)=S^{ab}\,\p_bf\p_ag(-1)^{\at\ft}
\end{equation}
where $S=\frac{1}{2}\,S^{ab}p_bp_a$ is an odd function on $T^*M$.
Parentheses stand for the canonical Poisson bracket on $T^*M$.
Clearly, see, e.g., \cite{tv:laplace1},   the Jacobi identity
for~\eqref{eqoddbracket} is equivalent to the equation
\begin{equation}\label{eqss}
    (S,S)=0
\end{equation}
on $T^*M$. (Notice that for an even symmetric bracket with the
even master Hamiltonian this would be an empty condition.)

Let us apply this to a long bracket on $M$, which is a usual bracket
on the manifold $\hat M$. Consider for simplicity long brackets of
weight $0$.

For an odd long bracket   specified by a master Hamiltonian $\hat\S$
by formula~\eqref{eqabstractlongbracket} (with $\l=0$) we see that it
satisfies the Jacobi identity if and only if $\hat\S$ satisfies the
equation $(\hat\S,\hat\S)=0$, with the canonical Poisson bracket on
$T^*\hat M$.

To get it more explicitly, we can express the canonical bracket on
$T^*\hat M$ in a ``$(D+1)$-formalism'', separating variables related
to $M$ from the extra variables $t,p_t$. For Hamiltonians on the
extended space $\hat M$ we have:
\begin{equation*}
(F,G)_{\hat M}=(F,G)_M+\der{F}{p_t}\der{G}{t}-\der{F}{t}\der{G}{p_t},
\end{equation*}
where for clarity we denoted by $(F,G)_{\hat M}$ the canonical
bracket on $T^*\hat M$ and by $(F,G)_{M}$ the  bracket on $T^*M$ with
$t,p_t$ considered as parameters. Applying this to
\begin{equation}\label{eqmasterham2}
    \hat \S=\frac{1}{2}\left(\S^{ab}p_bp_a+ 2t\g^ap_ap_t+t^2\lt
    p_t^2\right),
\end{equation}
we obtain the following theorem. Let
\begin{equation}\label{eqmasteronm}
    \S=\frac{1}{2} \,\S^{ab}p_bp_a
\end{equation}
be the master Hamiltonian for the bracket on $M$, and
\begin{equation}\label{eqgammaasham}
    \g=\g^{a}p_a
\end{equation}
be the local Hamiltonian specifying the upper connection in $\Vol M$.

\begin{thm}\label{thmjacobilong}
The Jacobi identity for an odd long bracket of weight zero
specified by the ``extended'' master
Hamiltonian~\eqref{eqmasterham2} is equivalent to the following
equations:
\begin{align}
    (S,S)=0 \label{eqssonM}\\
    (S,\g)=0 \label{eqswithgamma}\\
    (S,\lt)+(\g,\g)=0 \label{eqswiththeta}\\
    (\g,\lt)=0 \label{eqgammawiththeta}
\end{align}
\end{thm}

Let us analyze the geometric  meaning of these equations.

Equation~\eqref{eqssonM} means that the bracket on $M$ satisfies the
Jacobi identity, i.e., $M$ is an odd Poisson (Schouten) manifold with
the Poisson tensor $S$. It follows that the operator $D:=(S,\ )$ on
functions on $T^*M$ is a differential. Hence we can
rewrite~\eqref{eqswithgamma} as
\begin{equation}\label{eqdgamma}
    D\g=0,
\end{equation}
and~\eqref{eqswiththeta} as
\begin{equation}\label{eqDtheta}
    D\lt+(\g,\g)=0.
\end{equation}
Equation~\eqref{eqdgamma} is nothing but the condition of
\textit{flatness} for the upper connection $\g$. Let us explain.

For any upper connection associated with  an odd bracket,
curvature makes sense if the bracket is a Poisson bracket. Since
for ``curvature'' we need a Lie bracket (a commutator)  we cannot
discuss curvature  for upper connections associated with arbitrary
symmetric tensors $S^{ab}$, for example, when $S^{ab}$ is even.

\begin{rem}
The curvature of a usual (``lower'') connection is defined  either
in terms of the exterior differential of a local connection form
or by comparing the commutator of covariant derivatives with the
commutator of vector fields. For an upper connection or a
``contravariant derivative'' $\nabla^{\o}$,  curvature can be
introduced in terms of the bracket of $1$-forms that comes from a
Poisson structure.
\end{rem}
\begin{rem}\label{remkoszulbrack}
This bracket of $1$-forms is a particular case of the bracket on
arbitrary forms known as the ``Koszul bracket'' for the even
Poisson case. Any even Poisson structure induces an odd Koszul
bracket of forms such that  raising indices by the Poisson tensor
takes it to the canonical odd Schouten bracket of multivector
fields. The usual exterior differential $d$ is mapped to the
``Lichnerowicz operator'', i.e.  the Schouten bracket with the
Poisson tensor. Hence, for an even Poisson bracket, the
``curvature form'' of any associated upper connection can be
defines as a $2$-vector field. An analog of this construction can
be carried over for odd Poisson brackets. Such bracket induces an
even bracket on forms (an analog of the Koszul bracket), such that
it is mapped by raising indices by $S^{ab}$ to the canonical
Poisson bracket of functions on $T^*M$. In this case the exterior
differential $d$ is mapped to the operator $D=(S,\ )$. In
particular, both upper ``connection form'' and ``curvature form''
are Hamiltonians  linear and quadratic in $p_a$, respectively.
\end{rem}

\begin{ex}\label{exoddgammafromvolume}
Since $D^2=0$, the condition of flatness~\eqref{eqdgamma} will be
identically satisfied if $\g=-D\log\rho=-S^{ab}\,\p_b\log \rho$
for some $\rho$. This is exactly the case of the upper connection
coming from a volume form. Substituting this into~\eqref{eqDtheta}
we get $D\lt+(D\log\rho,D\log\rho)=0$, or
$D\bigl(\lt-(D\log\rho,\log\rho)\bigr)=0$, since $D$ is a
derivation of the canonical bracket. As the last term is
$\{\log\rho,\log\rho\}$, we have $\lt=\{\log\rho,\log\rho\}\ +\
(\text{\footnotesize Casimir functions})$. Notice
that~\eqref{eqgammawiththeta} is then satisfied identically.
\end{ex}

There is an important application to the case when the bracket of
functions on $M$ is non-degenerate. Then $M$ is an odd symplectic
manifold with the symplectic form
\begin{equation}\label{eqsymplform}
    \o=\frac{1}{2}\,dx^a\,dx^b\,\o_{ba}
\end{equation}
where $(\o_{ab})$ is the inverse matrix for $(S^{ab})$.  To the upper
connection $\g^a$ corresponds the usual connection
$\g_a=\o_{ab}\g^b$. The equation~\eqref{eqdgamma} for $\g^a$ is
equivalent to ${d(dx^a\g_a)=0}$, i.e., to the flatness of $\g_a$ in
the usual sense. Any flat connection in $\Vol M$ comes from a volume
form up to a  twist (see Remark~\ref{remtwistedvolumeform}), so we
can rewrite $\g_a$ as $\g_a=-\p_a\log \rho$. We find ourselves in the
situation of Example~\ref{exoddgammafromvolume}. In particular, for
the ``Brans--Dicke field'' $\lt$ we get
$\lt=\{\log\rho,\log\rho\}=\g^a\g_a$ (if we assume that there are no
odd constants). We arrive at the following theorem:

\begin{thm}[``Existence of action'']\label{thmexistenceofaction}
If $S^{ab}$ is non-degenerate, then the Jacobi identity for the long
bracket   is equivalent to the following conditions:\\
(1) $M$ is symplectic \\
(2) $\g^a=-S^{ab}\p_b\log\rho$ for some $\rho=e^{\Act}$ (a local
volume
form)\\
(3) $\lt=\g^a\g_a$.\\
\end{thm}

The logarithm of a volume form has the physical meaning of
``action''. Theorem~\ref{thmexistenceofaction} tells that if an odd
symplectic bracket can be extended to a long bracket of densities
satisfying the Jacobi identity,  there is an action function $\Act$
(at least, local) such that the long bracket comes from the volume
form $\bro=e^{\Act}\,Dx$ determined by this action.

For an arbitrary long bracket the coefficients  $S^{ab}$, $\g^a$
and $\lt$ are independent degrees of freedom. We see that the
Jacobi identity eliminates some degrees of freedom, reducing a
long odd Poisson bracket in the non-degenerate case to the Poisson
bracket of functions plus a volume form (in general, twisted) as
an extra piece of data.

The long bracket arising from a volume form $\bro$ is  simply
\begin{equation}\label{eqlongbrackfromvolume}
    \{\ps ;\ch \}=
    \{\bro^{-w_1}\ps,\bro^{-w_2}\ch\}\bro^{w_1+w_2}
\end{equation}
for  densities of weights $w_1$ and $w_2$, where at the r.h.s. stands
the bracket of functions defined by $S^{ab}$. Notice that the r.h.s.
of~\eqref{eqlongbrackfromvolume} is well-defined even for  local
volume forms with a non-trivial twist (see
Remark~\ref{remtwistedvolumeform}).

\subsection{$\D^2$ and flatness}

Let us first discuss  a purely algebraic situation.

Let $\D$  be a second order operator in a commutative associative
algebra $A$ with a unit. Consider the operator $\D^2$. In general,
$\ord\,\D^2\leq 4$. If, however, $\D$ is odd, then
$\D^2=\frac{1}{2}\,[\D,\D]$, hence $\ord \D^2\leq 3$. In the sequel
we   consider   odd operators. What is the meaning of the conditions
$\ord\,\D^2\leq k$ (for $k=2,1,0$)?

\begin{prop}\label{propabstorder2}
The condition $\ord\, \D^2\leq 2$ is equivalent to the Jacobi
identity for the   bracket generated by $\D$.
\end{prop}

Recall that the formula for the bracket  is (for odd $\D$)
\begin{equation}\label{eqoddbracket1}
     \{a,b\}=\D(ab)-(\D a)\,b-(-1)^{\at}a\,(\D b) +\D(1)\,ab.
\end{equation}
By redefining $\D$ as $\D-\D (1)$ the last term can be eliminated.

Proposition~\ref{propabstorder2} must be evident, since the
bracket~\eqref{eqoddbracket1} is nothing but the polarized principal
symbol of $\D$, and vanishing of $[\D,\D]$ modulo operators of order
$\leq 2$ is equivalent to the vanishing of the canonical Poisson
bracket of the principal symbol of $\D$ with itself.

\begin{prop}\label{propabstorder1}
The condition $\ord\, \D^2\leq 1$ is equivalent (in addition to the
Jacobi identity for the bracket) to the derivation property
\begin{equation}\label{eqderofbracket}
    \D \{a,b\}=\{\D a, b\}+(-1)^{\at}\{a,\D b\}.
\end{equation}
\end{prop}

Notice that if $\D(1)=0$, then $\ord\, \D^2\leq 1$ means that $\D^2$
is a derivation of the algebra $A$. Proposition~\ref{propabstorder1}
then means that $\D^2$ is a derivation of the associative
multiplication if and only if $\D$ is a derivation of the bracket.
(Then $\D^2$ is also a derivation of the bracket.)
See~\cite{yvette:divergence}.

Finally, if $\ord\,\D^2= 0$, %is of order $0$,
this basically means that $\D^2=0$. All previous properties hold,
i.e., $\D$ generates an odd Poisson bracket  for which it is a
derivation; and in addition $\D$ is a differential.

Odd Poisson algebras endowed with an operator $\D$ generating the
bracket and satisfying $\D^2=0$ have received the name of
\textit{Batalin--Vilkovisky algebras}. We see that imposing the
conditions $\ord\,\D^2\leq k$,   $k=2,1,0$,  allows to recover the
defining identities of the Batalin--Vilkovisky algebras step by
step.

Now let us return to the differential-geometric situation. Let $A$ be
the algebra of smooth functions on a manifold $M$. Any odd operator
$\D$ in $A$ of order $\leq 2$,
\begin{equation}\label{eqoddoperator}
    \D=\frac{1}{2}\,S^{ab}\,\p_b\p_a +T^a\,\p_a,
\end{equation}
is specified by a quadratic Hamiltonian $S=\frac{1}{2}\,S^{ab}p_bp_a$
and an associated upper connection $\g^a=\p_bS^{ba}-2T^a$ in the
bundle $\Vol M$ (compare~\eqref{eqdeltaonfunct}). (We have set for
convenience  $\D(1)=0$.)

Suppose  $\ord\D^2\leq 2$. That means $(S,S)=0$, and the odd bracket
generated by $\D$ makes $M$ into a Schouten ( = odd Poisson)
manifold.

By Proposition~\ref{propabstorder1}, $\D$ is a derivation of the
Schouten bracket if and only if $\D^2$ is a vector field (which is,
moreover, a Poisson vector field).

\begin{thm}\label{thmflatnessupper}
Let $\D$ be an arbitrary odd second order
operator~\eqref{eqoddoperator}, such that $\D(1)=0$. The properties
that $\D$ is a derivation of the corresponding Schouten bracket and
that $\D^2$ is a Poisson vector field, are equivalent to the flatness
of the upper connection:
\begin{equation}\label{eqflatness}
    D\g=0,
\end{equation}
where $D=(S,\ )$ and  $\g=\g^ap_a$.
\end{thm}

\begin{cor}\label{corflatgn}
Let the upper connection $\g^a$ come from a usual connection $\g_a$
in $\Vol M$ as $\g^a=S^{ab}\g_b$.  Then  $\D$ is a derivation of the
Schouten bracket (equivalently, $\D^2$ is a Poisson vector field) if
and only if $\g_a$ is flat  ``in the directions of Hamiltonian vector
fields'', i.e. $S^*(d\gn)=0$ where $S^*(dx^a)=S^{ab}p_b$ and
$\gn=\g_a\,dx^a$. In particular, $\D$ is a derivation of the Schouten
bracket  in the case if  $\g_a$ is flat.
\end{cor}

Indeed,  $D\g=S^*(d\gn)$, as raising indices by $S^{ab}$ takes $d$ to
$D$, see Remark~\ref{remkoszulbrack}.

\begin{rem}
A statement equivalent to Corollary~\ref{corflatgn} was obtained  in
the important paper~\cite{yvette:divergence}.  Notice that the class
of operators~\eqref{eqoddoperator} where $\g^a$ comes from a
connection in $\Vol M$ coincides with the class of operators given by
``abstract divergences''  (see
Remark~\ref{remlapwithgenuineconnect}). Divergence operators in
$\Vect (M)$ and connections in $\Vol M$ are equivalent notions. In
particular, a connection is flat if the corresponding divergence is
flat, i.e., satisfies~\eqref{eqdivkommutator}.
\end{rem}

Consider now the algebra of densities $\V(M)$ endowed with the
canonical scalar product~\eqref{eqscprod2} and an odd operator $\D$
in it. Suppose $\D^*=\D$ and $\D(1)=0$. Such operators (in other
words, odd canonical pencils) are in $1-1$ correspondence with odd
brackets in $\V(M)$, i.e., odd long brackets on $M$
(Theorem~\ref{thmmain}).

\begin{thm} \label{thmdeltakvadratpencil}
Let $\D_w$ be the canonical pencil corresponding to an odd long
bracket on $M$. If the long bracket satisfies the Jacobi identity,
then
\begin{equation*}%\label{}
    \D_w^2=\lie{X}
\end{equation*}
where $X$ is a Poisson vector field on $M$. Here $\lie{X}$ stands for
the Lie derivative on $w$-densities.
\end{thm}

(Let us emphasize that $X$ is a vector field on $M$, which is Poisson
w.r.t. the odd bracket of functions.)

Indeed, by Proposition~\ref{propabstorder2}, the Jacobi identity for
a long bracket is equivalent to $\ord\D^2\leq 2$. However, since $\D$
can be constructed as the Laplace operator  corresponding to the
canonical divergence on $\hat M$ (see Theorem~\ref{thmlaplacian}),
there is a ``jump'':  if $\ord\,\D^2\leq 2$, then automatically
$\ord\,\D^2\leq 1$. This is an algebraic fact following from the
flatness of the canonical divergence~\eqref{eqdivkommutator}; compare
with Corollary~\ref{corflatgn} and the remark after it. Hence
$\D^2=\X$ is a vector field on $\hat M$.

\begin{lm} The vector field $\X=\D^2$ on $\hat M$ is divergence-free:
\begin{equation}\label{eqdivfree}
    \div \X=0,
\end{equation}
where $\div$ in~\eqref{eqdivfree} is the canonical
divergence~\eqref{eqdivfordensities}.
\end{lm}

Not proving this simple lemma (compare Proposition 2.2
in~\cite{tv:laplace1}), we conclude, by
Theorem~\ref{thmclassderonhatm}, that $\X=\lie{X}$ for a vector field
on $M$ (which must be Poisson), and
Theorem~\ref{thmdeltakvadratpencil} follows.

\begin{ex}
Consider the symplectic case, i.e., when the bracket on $M$ is
non-degenerate. Then, by Theorem~\ref{thmexistenceofaction}, the
Jacobi identity for the long bracket implies that $\g^a$ comes from a
flat connection $\g_a=-\p_a\log\rho$ in $\Vol M$, and the canonical
pencil $\D_w=\D^{LB}_w$ is the Laplace--Beltrami pencil of
Example~\ref{exLB}. Then $\D_w^2=\lie{X}$ for the Hamiltonian vector
field
\begin{equation}\label{eqxcanon}
    X=\grad \,\frac{\D_{\text{can}}(\bro^{1/2})}{\bro^{1/2}},
\end{equation}
where $\D_{\text{can}}$ stands for the canonical Laplacian on
half-densities~\cite{hov:semi, hov:proclms}. Notice that the
Hamiltonian in~\eqref{eqxcanon} is well-defined even if the volume
form exist only locally. A condition $\D_w^2=0$ will be equivalent to
the Batalin--Vilkovisky equation
\begin{equation}\label{eqBV}
    \D_{\text{can}}(\bro^{1/2})=0
\end{equation}
or
\begin{equation}\label{eqBV2}
    \D_{\text{can}}\,\left(e^{\Act}\,Dx\right)^{{1}/{2}}=0
\end{equation}
if $\bro=e^{\Act}$.
\end{ex}

\section{Generalizations}\label{secgen}

\subsection{Operators and brackets of non-zero weight}

The study of operators and brackets of non-zero weight can have a
considerable interest.

Let us briefly indicate some of the statements and formulae for
$\l\neq 0$ which we have omitted in the previous sections.

For a canonical pencil $\D_w$ the singular points will be $w=0$,
$w=\frac{1-\l}{2}$, $w=1-\l$. In particular, the role of
half-densities  will be played by densities of weight
$\frac{1-\l}{2}$ for general $\l$.

$S^{ab}$ is no longer a tensor field; it is a tensor density of
weight $\l$. Similar is true for $\g^a$. It is an upper
``connection-density''.  In contrast,   $\g_a$ (when it appears) has
to be a genuine connection; it will acquire weight by raising indices
with the help of $S^{ab}$.

An interesting example of operators of non-zero weight acting on
densities is well known in integrable systems. Namely, consider a
one-dimensional manifold $M$ with a coordinate $x$, and set $\l=2$.
By inspection of formulae~\eqref{eqzamenaSl}--\eqref{eqzamenathetal}
we see that the parameters $s,\g,\lt$ of the canonical pencil  (which
we have for convenience  multiplied by $2$)
\begin{equation*}%\label{}
    \D_w=s\,\p^2+ \bigl(s_x+(2 w+1)\g\bigr)\,\p + w\,\g_x  + w(w+1)\,\lt,
\end{equation*}
where $\p=d/dx$, under a change of a coordinate $y=y(x)$  transform
as follows:
\begin{align}
    s'&=s
    \label{eqzamenas}
    \\
     \g'&=\frac{y_x\g+y_{xx}s}{(y_x)^2}
    \label{eqzamenag}
    \\
    \lt'&=\frac{(y_x)^2\lt+2y_xy_{xx}\g+(y_{xx})^2s}{(y_x)^4}
    \label{eqzamenath}
\end{align}
We used  dashes for the parameters in a new coordinate system. As $s$
is invariant, we can set $s=1$. Then we get
\begin{equation*}%\label{}
    \D_w=\p^2+  (2 w+1)\,\g \,\p + w\,\g_x + w(w+1)\,\lt
\end{equation*}
where $\g_x=d\g/dx$. In particular,  for $w=-\frac{1}{2}$ we obtain a
`Sturm--Liouville' operator
\begin{equation}\label{eqoursturm}
    L=\D_{-\frac{1}{2}}=\p^2    -\frac{1}{2}\, \left(\g_x +\frac{1}{2}\,\lt\right)
\end{equation}
with the potential $U=\frac{1}{2}\, \left(\g_x
+\frac{1}{2}\,\lt\right)$. It maps densities of weight $-\frac{1}{2}$
to densities of weight $\frac{3}{2}$.  A known fact about such
Sturm--Liouville operators $L=\p^2-U$ is that $U$  has a
transformation law involving the Schwarz derivative:
\begin{equation*}
    U'=(y_x)^{-2}\left(U-\frac{1}{2}\,\Schw[y(x)]\right),
\end{equation*}
where
$\Schw[y(x)]=\frac{y_{xxx}}{y_x}-\frac{3}{2}\,\bigl(\frac{y_{xx}}{y_x}\bigr)^2$
(see, e.g., \cite{hitchin:hsw}). In our
parametrization~\eqref{eqoursturm}, it   follows from the
transformation laws~\eqref{eqzamenag}, \eqref{eqzamenath} for $\g$
and $\lt$, which are a very special case of the transformation law
for the coefficients of a long bracket.

It will be interesting to study the geometrical meaning of this
relation further, in particular, to explore  possible links of our
constructions with projective connections.

Speaking about odd brackets of weight $\l\neq 0$, we can notice
that in the odd case ``Poisson brackets of functions'' make no
sense separately from the brackets of densities, since the bracket
of functions is a density of weight $\l$, so the Jacobi identity
has to involve densities of all weights.
Theorem~\ref{thmjacobilong} is generalized for arbitrary $\l$ as
follows:

\begin{thm}\label{thmjacobilonglambda}
The Jacobi identity for an odd long bracket of weight $\l$
specified by the master Hamiltonian~\eqref{eqmasterham} is
equivalent to the following equations:
\begin{align}
    (S,S)=2\l S\g \label{eqssonMl}\\
    (S,\g)=\l S\lt \label{eqswithgammal}\\
    (S,\lt)+(\g,\g)=\l \g \lt \label{eqswiththetal}\\
    (\g,\lt)=0 \label{eqgammawiththetal}
\end{align}
\end{thm}

In  equations~\eqref{eqssonMl}--\eqref{eqgammawiththetal} the
quantities  $S$, $\g$ and $\lt$ are considered as functions on $T^*M$
whose definition depends on coordinates or on a choice of volume form
on $M$.  In particular, as for $\l\neq 0$  the Hamiltonian  $S$ takes
values in $\l$-densities, the canonical bracket $(\S,\S)$ does not
have an invariant meaning. This is indicated by the presence of a
non-invariant term $S\g$ in the r.h.s. of~\eqref{eqssonMl}.

It will be interesting to study for $\l\neq 0$ the examples of a
non-degenerate $S^{ab}$  (a ``symplectic structure taking values in
densities'') and a long bracket coming from a connection or a volume
form.

\subsection{Operators of higher order}
\label{subsechigh}

It is tempting to extend the classification of operators and brackets
that we have obtained for the operators of order $\leq 2$ to
operators of higher order. At the moment, we have more questions than
answers concerning this case.

First of all, the algebraic framework for brackets generated by an
operator $\D$ is as follows. For simplicity of notation let $\D$
be even and $A$ be purely even. Of course, the interesting case is
that of an odd $\D$. Recall that an operator $\D$ acting in a
commutative associative algebra $A$ has order $\leq n$ if and only
if all $(n+1)$-fold commutators
$[\ldots[[\D,a_1],a_2],\ldots,a_{n+1}]$ vanish (where $a_i$ are
arbitrary elements of $A$). Define a sequence of  ``higher
brackets'' corresponding to the operator $\D$ as $k$-fold
commutators with elements of $A$ applied to $1$:
\begin{equation*}
    \begin{split}
    \{a\}&=[\D,a](1)=\D a-\D(1) \cdot a   \\
    \{a,b\}&=[[\D,a],b](1)=\D(ab)-\D(a)\cdot b-a\cdot \D(b)+\D(1)\cdot
    ab   \\
    \{a,b,c\}&=[[[\D,a],b],c](1)=
    \D(abc)-\D(ab)\cdot c
    -\D(ac)\cdot b-\D(bc)\cdot a\\
     \quad &\quad  +
    \D a\cdot bc+\D b\cdot ac+\D c\cdot ab-\D(1)\cdot
    abc\\
    \ldots
    \end{split}
\end{equation*}
up to the $n$-fold bracket, which simply coincides with the
$n$-fold commutator. The brackets higher than the $n$-th vanish.
Notice that all the brackets $\{a_1,\ldots,a_k\}$ are symmetric.
Each of them  is the obstruction for the previous bracket to be a
derivation (w.r.t. each of the arguments). The $k$-th bracket is
an operator of order $n-k+1$ on each of its arguments. The top
$n$-th bracket $\{a_1,\ldots,a_n\}$ is a multi-derivation. It is
exactly the polarized principal symbol of $\D$, i.e., the
principal symbol considered as a symmetric multilinear function.
We can call the constructed sequence of higher brackets the
sequence of \textit{polarizations} of the operator $\D$.

For example, for an operator of order $\leq 2$ we essentially have
just one bracket,
\begin{equation*}
    \{a,b\} =[[\D,a],b](1)=\D(ab)-\D(a)\cdot b-a\cdot \D(b)+\D(1)\cdot
    ab,
\end{equation*}
which is a bi-derivation. It was our main  object  in the previous
sections. (The unary  bracket  $\{a\}=[\D,a](1)=\D a-\D(1)\cdot a$
was used when we redefined $\D$ so to have $\D(1)=0$.)

\begin{rem}
The definition of the order of an operator in terms of commutators is
traced back to Grothendieck. The above sequence of polarizations for
$\D$ essentially coincides with the operations $\Phi_k$ considered by
Koszul in~\cite{koszul:crochet85}.
\end{rem}

For an odd operator $\D$ of order $\leq n$, imposing the
conditions $\ord \D^2\leq r$, with $r=2n-2, \,2n-3, \ldots$, as in
Section~\ref{secjacobi}, will lead to an hierarchy of identities
involving the associative multiplication, brackets of various
orders and $\D$. The resulting structure  can be loosely called a
``homotopy Batalin--Vilkovisky algebra''. This should be further
explored.

It is instructive to write down the higher brackets explicitly for
an operator of order $\leq 3$ in a differential-geometric setting.
Consider
\begin{equation}\label{eqoperthird}
    \D=\frac{1}{6}\,S^{abc}\,\p_c\p_b\p_a+
    \frac{1}{2}\,P^{ab}\,\p_b\p_a+T^a\,\p_a+R.
\end{equation}
Let us set $R=0$ for simplicity. Then $\{f\}=\D f$,
\begin{equation*}%\label{}
    \{f,g\}=\frac{1}{2}\,S^{abc}\bigl(\p_cf\cdot \p_b\p_ag+
    \p_c\p_bf\cdot \p_ag\bigr) +P^{ab}\,\p_bf\cdot \p_ag
\end{equation*}
and
\begin{equation*}%\label{}
    \{f,g,h\}=S^{abc}\,\p_cf\,\p_bg\,\p_ah.
\end{equation*}
An analog of the question considered in the previous sections is
how to recover $\D$ from the binary and ternary brackets,
$\{f,g\}$ and $\{f,g,h\}$. Clearly, the problem is to identify in
geometrical terms extra data that can be constructed from $\D$ and
together with the brackets will allow to recover $\D$.

For operators of the second order such information is contained in
the subprincipal symbol, which we are able to interpret as an upper
connection. However, for the an operator of the third
order~\eqref{eqoperthird}, this will not give a solution, since $T^a$
will not enter $\sub \D$. One could look for something like iterated
``sub-subprincipal'' symbols, possibly defined using several volume
forms.

%In general, the situation for operators of the even and odd order
%might differ, and these cases should be considered separately. We
%hope to study this problem elsewhere.

With this might be related the following construction for operators
of the second order: consider the subprincipal symbol as a vector
field depending on a connection and take its divergence w.r.t. the
other connection. Notice that for operators defined by a (genuine)
connection this is an object depending on three connections. Taking
one of them at the midpoint of the segment joining the other two, we
arrive at  a construction of a scalar function from two genuine
connections in $\Vol M$ if a second order principal symbol $S$ is
given:
\begin{equation}\label{eqdivv}
   \coc{0}{1}= \divv{0}{1}=\p_a(\g_1^a-\g_0^a)-
    \frac{1}{2}\,(\g_{0a}+\g_{1a})(\g_1^a-\g_0^a)
\end{equation}
(indices are raised  by the tensor $S^{ab}$). It is exactly the
`Batalin--Vilkovisky' term $\div_{\g}X-\frac{1}{2}\,X_aX^a$ appearing
in the equation~\eqref{eqchangeofdeltaunderconnection} for the
transformation of the canonical pencil corresponding to a connection
under a change of connection. (Notice that a convex combination of
connections like $\frac{1}{2}\,(\g_{0}+\g_{1})$ makes good sense.)

The function~\eqref{eqdivv} satisfies the cocycle property
\begin{equation*}%\label{}
    \coc{0}{1}+\coc{1}{2}=\coc{0}{2}
\end{equation*}
or
\begin{equation*}%\label{eqstrannaja}
    \divv{0}{1}+\divv{1}{2}=\divv{0}{2}
\end{equation*}
which is just  the groupoid property of the `Batalin--Vilkovisky
equations'~\eqref{eqbv1}.

\bigskip \textbf{Acknowledgements:}
Cordial thanks  go to Yvette Kosmann-Schwarzbach, Kirill Mackenzie
and M.~A.~Shubin for numerous stimulating discussions.

\def\cprime{$'$}

\end{document}